\documentclass[10pt]{article}

\usepackage{amsfonts}
\usepackage{amssymb}
\usepackage{amsmath}
\usepackage{latexsym}
\usepackage{graphicx}
\usepackage{amsmath,bm}

\usepackage{mathrsfs,amsfonts,dsfont}
\usepackage{hyperref}
\usepackage{color}
\usepackage{geometry}
\usepackage{multicol}
\usepackage{multirow}
\usepackage{colortbl}
\usepackage[bottom]{footmisc}
\usepackage{url}
\usepackage{datetime}

\usepackage{bm}
\usepackage{jabbrv}

\usepackage[linesnumbered,ruled,vlined]{algorithm2e}
\usepackage{algorithmicx}
\usepackage{algpseudocode}
\usepackage{varwidth}

\newcommand{\V}[1]{\boldsymbol{#1}}

\newtheorem{theorem}{Theorem}[section]

\newcommand{\be}{\begin{equation}}
\newcommand{\ee}{\end{equation}}
\newcommand{\ba}{\begin{array}}
\newcommand{\ea}{\end{array}}
\newcommand{\ben}{\begin{eqnarray}}
\newcommand{\een}{\end{eqnarray}}
\newcommand{\bn}{\begin{eqnarray*}}
\newcommand{\en}{\end{eqnarray*}}

\def\epsilonb {\boldsymbol{\epsilon}}

\begin{document}
\centerline{\Large \bf Comparisons of Some Iterative
Algorithms for Biot Equations}

\renewcommand{\thefootnote}{\fnsymbol{footnote}}

\vskip .5cm

\centerline{
Dedicated to the Memory of Professor V. V. Zhikov
}

\vskip .5cm
\renewcommand{\thefootnote}{\arabic{footnote}}
\centerline{
{\bf Mingchao Cai}$\ddag$\footnote{Corresponding author. E-mail address:
cmchao2005@gmail.com. This author's work is supported in part by NIH BUILD grants through ASCEND Pilot project and NSF HBCU-UP Research Initiation Award through HRD1700328.}
{\bf Guoping Zhang}$\ddag$\footnote{E-mail address:
guoping.zhang@morgan.edu. This author's work is supported in part by NIH BUILD grants through ASCEND Pilot project.}
}

\vskip .3cm \centerline{$\ddag$\it Department of Mathematics, Morgan State University,}
\centerline{\it 1700 E Cold Spring Ln, Baltimore, MD 21251, USA.}

\vskip .5cm \noindent\rule[2mm]{\textwidth}{.1pt}

\noindent {\bf Abstract}
\vskip .5cm

In this paper, we aim at solving the Biot model under stabilized finite element discretizations.
To solve the resulting generalized saddle point linear systems, some iterative methods are proposed and compared. In the first method, we apply the GMRES algorithm as the outer iteration. In the second method, the Uzawa method with variable relaxation parameters
is employed as the outer iteration method. In the third approach, Uzawa method is treated as a fixed-point iteration, the outer
solver is the so-called Anderson acceleration. In all these methods, the inner solvers are preconditioners  for the generalized saddle point problem. In the preconditioners, the Schur complement approximation is derived by using Fourier analysis approach. These preconditioners are implemented exactly or inexactly. Extensive experiments are given to justify the performance of the proposed preconditioners and to compare all the algorithms.

\vskip .5cm \noindent
{AMS: 65M55, 65N22, 65N30, 65F08, 65M22}

\noindent\rule[2mm]{\textwidth}{.1pt}

\section{Introduction}
Biot's equations, which describe the deformation of poroelastic material under saturated fluids, have many wide applications in geosciences
and biomechanics. For example, Biot model is frequently used in simulating brain swelling model to quantify brain edema assessment. Combined with image data and patient-specific data such as cerebral blood flow conditions, numerical methods are applied for simulating brain swelling under ischemic conditions or after traumatic brain injury. In practical applications, numerical methods such as Finite Element methods and Finite Difference methods are normally used in simulations. As the model contains several physical parameters, numerical study for this model is very
challenging. To capture the correct behavior of the numerical solution, the discretization schemes usually require extremely fine grids which
will lead to large-scale generalized saddle point discrete linear systems. This paper aims at providing efficient numerical methods for solving the resulting linear system of the Biot model under stabilized Finite Element discretizations.

For simplicity, we assume that the domain $\Omega$ is a unit square or a unit cube so that we can avoid the rescaling of the spatial variables.
The (quasi-static) Biot equations is
\begin{equation}\label{Biot_model}
 \left \{
 \begin{array} {ll}
- \mbox{div} (2\mu \epsilonb(\bm{u})) - \nabla \lambda \mbox{div} \bm{u}  +  \alpha \nabla p= \bm{f}, \quad & \forall \bm{x} \in \Omega,\\
\frac{\partial}{\partial t} (s p+  \alpha \mbox{div} \bm{u}) -\mbox{div} (\kappa \nabla p ) = g, \quad & \forall \bm{x} \in \Omega.
\end{array}
    \right.
\end{equation}
In the above equations, $\bm{u}$ is the displacement of the elastic deformable medium,
$$
\epsilonb(\bm{u}) = \frac{1}{2} \left[\nabla \bm{u} + \nabla^T \bm{u} \right],
$$
$p$ is the pressure of fluid, $\bm{f}$ is the body force, $g$ is a source term for fluid, $\lambda$ and $\mu$ are Lam\'{e} constants, which can be expressed as
\begin{equation}\label{Poisson_exp}
\lambda = \frac{E \nu}
{(1+\nu)(1-2\nu)},
\quad  \mu=  \frac{E}
{2(1+\nu)},
\end{equation}
$s \ge 0$, is the constrained specific storage coefficient, $\kappa >0$ represents the permeability, $\alpha$ is the Biot constant which is
 close to 1. In (\ref{Poisson_exp}), $E$ is the Young's module, $\nu$ is the Poisson ratio. When $\nu$ approaches to $0.5$, the elastic material is almost incompressible. Note that $\alpha$ is close to $1$ and one can rescale $\alpha p = \tilde{p}$, for simplicity and without loss of generality, we assume $\alpha=1$ so that the final system is symmetric.
For the boundary condition, one can apply different types of boundary conditions. For simplicity, we will apply the typical homogeneous Dirchlet boundary condition for both the displacement and pressure \cite{gaspar2003finite, gaspar2007efficient}. The initial condition should satisfy the Stokes equations \cite{lee2017parameter, murad1996asymptotic, feng2017analysis}. Except specifically pointing out, we will focus on the discussions on pure Dirichlet boundary conditions. Other type of boundary conditions like mixed Dirichlet and Neumann boundary condition \cite{murad1996asymptotic} will also be discussed when necessary.

 For the second equation of (\ref{Biot_model}), we apply the backward Euler scheme:
\begin{equation}\label{Euler_scheme}
- (\mbox{div} \bm{u}^{n+1} - \mbox{div} \bm{u}^{n} ) - s (p^{n+1}-p^n) + \Delta t \mbox{div} (\kappa \nabla p^{n+1}) = - \Delta t g^{n+1}.
\end{equation}
By the backward Euler scheme (\ref{Euler_scheme}), the resulting differential operator is
\begin{equation}\label{mixed_operator}
\mathcal{M} = \left[\begin{array}{ccc}
-\mbox{div} (2\mu \epsilonb(\cdot))   &  \mbox{grad}    \\
-  \mbox{div}   & -(s I - \theta  \Delta) \\
\end{array}
\right].
\end{equation}
Here, $\theta= \kappa \Delta t$.

For spatial discretization, we consider to use the stabilized Mini element \cite{rodrigo2016stability} or stabilized $P_1$-$P_1$ (or $Q_1-Q_1$) discretization. In the literature, Gaspar et. al. develop finite difference discretizations \cite{gaspar2003finite, gaspar2007efficient, naumovich2006finite}, Lipnikov uses $P^1$ elements for displacement and a lowest-order Raviart-Thomas elements for pressure \cite{lipnikov2002numerical}, Lee et. al apply the conforming and nonconforming finite elements \cite{lee2016robust, lee2017parameter}. 
We adopt the stabilized Finite Element method proposed in \cite{rodrigo2016stability} because this method leads to monotone scheme which is very important to suppress the pressure approximation errors. After the discretization, the resulting linear system is of the form:
\begin{equation}\label{stab_sys}
{\bf M} \left[\begin{array}{c}
{\bf u} \\
p
\end{array}\right]=
\left[\begin{array}{cc}
{\bf A}   & {\bf B}^t \\
{\bf B}   & -{\bf D}
\end{array}\right]
\left[\begin{array}{c}
{\bf u} \\
{\bf p}
\end{array}\right]=
\left[\begin{array}{c}
 {\bf f} \\
 {\bf g}
\end{array}\right].
\end{equation}
In designing fast solvers for such a generalized saddle point problem, there are two main ingredients: preconditioners and iterative methods. In the following, we highlight the contribution of this paper on both the preconditioning aspect and the iterative method aspect.

For preconditioning the Biot system, the main focus of this paper is the approximation of the Schur complement based on a Fourier analysis approach. In a recent paper \cite{lee2017parameter}, the authors consider to introduce an intermediate variable (one can call it "total pressure", which is a linear combination of $\mbox{div} \bm{u}$ and the fluid pressure), and reformulate the $2$-by-$2$ saddle point operator (\ref{mixed_operator}) into a $3$-by-$3$ saddle point problem, then they use the conforming Finite Element method to approximate the weak form in the functional space $(H^1_0(\Omega))^d \times  L^2(\Omega) \times  H^1_0(\Omega)$. Their arguments are based on that proper functional spaces can be naturally introduced which correspond to the new prime variables and uniform inf-sup stability conditions hold naturally true on both continuous and discrete levels. We note that introducing one more variable will increase the problem size and approximating the boundary condition for the intermediate variable is sophisticated (cf. \cite{lee2017parameter} for the details). In contrast, we keep the original $2$-by$2$ saddle point form and apply the stabilized finite element discretizations \cite{rodrigo2016stability}. Moreover, in investigating the robustness of the proposed method with respect to parameters, the number of parameters involved in is as less as possible.

For the iterative methods, we will consider three different approaches: a) a preconditioned GMRES method \cite{saad1986gmres, saad1993flexible} applied to the generalized saddle
point systems; b) a variable relaxation parameter Uzawa algorithm \cite{hu2001iterative}, in which the relaxation parameters are selected for minimizing the energy-norm errors of each substep; c) an Anderson acceleration algorithm \cite{walker2011anderson, toth2015convergence, ho2015accelerating}, applied to a fixed point formulation of
a Uzawa algorithm for the generalized saddle point system. In all these three algorithms, preconditioners are the same for both the $(1,1)$ block and the Schur complement. By carefully designing and optimizing the parameters for all of the key components of the solvers, ranging from the inexact solvers to the robustness with respect to
the physics parameters, we compare the solvers that can readily  be employed to design the novel algorithms for FE spatial discretizations of the equations of the Biot model.

The organization of the paper is as follows. In Section 2, we introduce the weak forms and the finite element approximations. In Section 3, preconditioners are introduced. The spectral properties of the linear operators and the preconditioned system are analyzed. In Section 4, the three iterative algorithms are presented. Extensive numerical experiments are given in Section 5 and conclusions are drawn in Section 6.

\section{Weak forms and finite element approximations}

For simplicity and without loss of generality, we focus on the constant-coefficient case. The functional spaces for the displacement and the fluid pressure are given by
$$
\begin{array}{lll}
\bm{V} &:=& \{\bm{v} \in (H^1(\Omega))^d :~\bm{v}|_{\partial\Omega} = \bm{0}. \},  \\
Q &:=& H^1_0 (\Omega) := \{p \in H^1(\Omega), ~p|_{\partial \Omega}  =0.\}.
\end{array}
$$
Multiplying the first equation of the Biot's model by a test function $\bm{v} \in \bm{V}$ and the second equation by a test function $q \in Q$, we have following weak problem: find $\bm{u} \in \bm{V}, p \in Q$, such that
 \begin{equation} \label{Biot_weakform}
    \left \{
        \begin{array} {ll}
      a(\bm{u},\bm{v})+b(\bm{v}, p)=(\bm {f},\bm{v})\qquad & \forall\bm{v}\in\bm{V},\\
      b(\bm{u},q)-d(p, q)=(g, q) \qquad & \forall q\in Q.
        \end{array}
    \right.
 \end{equation}
Here, the bilinear forms are
\begin{equation}
a(\bm{u}, \bm{v})  =
\int_{\Omega} 2 \mu \epsilonb(\bm{u}):\epsilonb(\bm{v}) + \lambda {\rm div}\bm{u} {\rm div}\bm{v}  d \Omega,
\label{a_form1}
\end{equation}
with
$$
\displaystyle
\epsilonb(\bm{u}) :\epsilonb(\bm{v}) =
\sum_{i=1}^d
\sum_{j=1}^d
\epsilonb_{ij}(\bm{u}) \epsilonb_{ij}(\bm{v}),
$$
\begin{eqnarray}
d(p, q) & = &
 \int_{\Omega} s p q + \theta \nabla p \cdot \nabla q  d \Omega,
\label{bilinear_forms_c}\\
b(\bm{v}, q) & = &   -\int_{\Omega}{\rm div}\bm{v}\ q \  d \Omega .
\label{bilinear_forms_b}
\end{eqnarray}
The functional space $\bm{V}$ can be endowed with the conventional $H^1$- norm, or the norm induced by the bilinear form $a(\cdot, \cdot)$ \cite{xu1992iterative, mardal2011preconditioning}. The functional space $Q$ can be endowed with the $H^1$- norm or the norm induced by $d(\cdot, \cdot)$. If $\lambda$ and $\mu$ are constant, there holds the following identity.
\begin{equation}\label{grad_div_iden}
- \mbox{div} (\mu [\nabla \bm{u}+\nabla \bm{u}^T]) - \nabla \lambda \mbox{div} \bm{u}= - \mu {\bf \Delta} \bm{u}- (\mu+\lambda) \nabla \mbox{div} \bm{u}.
\end{equation}
Then, starting from the right hand side of (\ref{grad_div_iden}), we have
$$
a(\bm{u}, \bm{v})  =
\int_{\Omega} \mu \nabla \bm{u}: \nabla \bm{v}  + (\mu+\lambda) {\rm div}\bm{u} {\rm div}\bm{v}  d \Omega.
$$

In the following, we shall use the Poincare's inequality
$$
{||q||}_{L^2} \leq {C}_P {||\nabla q||}_{L^2},
$$
and the first Korn inequality
\begin{equation}\label{Korn_ineq}
{||\epsilonb(\bm{v})||}_{L^2} \geq {C}_K {||\bm{v}||}_{H^1} \quad  \forall \bm{v} \in \bm{V}.
\end{equation}
For $\forall \bm{v} \in \bm{V}$, there also hold
\begin{eqnarray}\label{div_est}
2 \mbox{div} \epsilonb(\bm{v}) - \nabla \mbox{div} \bm{v} = {\V \Delta} \bm{v} = \nabla \mbox{div} \bm{v} - \mbox{curl}^2 \bm{v}, \cr
|| \mbox{curl} \bm{v}||^2_{L^2} + ||\mbox{div} \bm{v}||^2_{L^2}  = ||\nabla \bm{v}||^2_{L^2} = 2 ||\epsilonb(\bm{v})||^2_{L^2} - ||\mbox{div} \bm{v}||^2_{L^2}, \cr
||\mbox{div} \bm{v}||_{L^2} \le ||\epsilonb(\bm{v})||^2_{L^2} \le ||\nabla \bm{v}||^2_{L^2}.
\end{eqnarray}
Noting from (\ref{Korn_ineq}) and (\ref{div_est}) that
$$
a(\bm{u} , \bm{u}) = 2\mu {{||\nabla \bm{u}||}^2}_{L^2} + \lambda {{||\mbox{div} \bm{u}||}^2}_{L^2} \geq  2 \mu {C}_K^2 {{||\bm{u}||}^2}_{H^1},
$$
we see that
\begin{equation}\label{a_coer}
|a(\bm{u},\bm{u})| \geq {\alpha}_A {{||\bm{u}||}^2}_{H^1}, \quad \mbox{with} \quad {\alpha}_A = 2\mu {C}_K^2.
\end{equation}
The bilinear form $a(\cdot, \cdot)$ is also bounded. More clearly,
\begin{equation}\label{a_bound}
\begin{array}{ll}
|a(\bm{u},\bm{v})| & \leq 2\mu {||\epsilonb (\bm{u})||}_{L^2} {||\epsilonb (\bm{v})||}_{L^2} + \lambda{||\mbox{div} \bm{u}||}_{L^2} {||\mbox{div} \bm{v}||}_{L^2}\\
& \leq 2\mu {||\bm{u}||}_{H^1} {||\bm{v}||}_{H^1} + \lambda {||\bm{u}||}_{H^1} {||\bm{v}||}_{H^1} \\
& \leq {C}_A {||\bm{u}||}_{H^1} {||\bm{v}||}_{H^1}
\end{array}
\end{equation}
where, $C_A = (2 \mu + \lambda)$. For the bilinear form $b(\cdot, \cdot)$, we have
\begin{equation}\label{b_bound}
\begin{array}{ll}
|b(\bm{u},q)| & \leq {||\mbox{div} \bm{u}||}_{L^2} {||q||}_{L^2} \leq {||\nabla \bm{u}||}_{L^2} {||q||}_{L^2}\\
& \leq {||\bm{u}||}_{H^1} {||q||}_{L^2} \leq C_P {||\bm{u}||}_{H^1} {||q||}_{H^1} \\
& \leq {C}_B {||\bm{u}||}_{H^1} {||q||}_{H^1}
\end{array}
\end{equation}
where, ${C}_B = C_P$. For the bilinear form $d(\cdot, \cdot)$, we have
\begin{equation}\label{d_bound}
\begin{array}{ll}
|d(p,q)| & \leq |s(p,q)| + |\theta(\nabla p, \nabla q)|   \\
& \leq s {||p||}_{L^2} {||q||}_{L^2} + \theta{||\nabla p||}_{L^2} {||\nabla q||}_{L^2} \\
& \leq s {C}_P {||\nabla p||}_{L^2} {||\nabla q||}_{L^2} + \theta{||\nabla p||}_{L^2} {||\nabla q||}_{L^2} \\
& \leq {C}_D {||p||}_{H^1} {||q||}_{H^1}
\end{array}
\end{equation}
with ${C}_D = s {C}_P + \theta$. On the other hand, we have
\begin{equation}\label{d_coer}
\begin{array}{ll}
d(p,p) & \geq s {||p||}_{L^2} + \theta {||\nabla p||^2}_{L^2} \\
& \geq \theta {||\nabla p||^2}_{L^2} \\
& \geq {\alpha}_D {||p||^2}_{H^1}
\end{array}
\end{equation}
with ${\alpha}_D = \theta {C}_K$.

In summary, from (\ref{a_coer}) to (\ref{d_coer}), for the generalized saddle point operator (\ref{mixed_operator}), the linear operators induced by $a(\cdot, \cdot),b(\cdot, \cdot),d(\cdot,\cdot)$ are bounded, and $a(\cdot, \cdot)$ and $d(\cdot, \cdot)$ are coercive. We apply the stabilized Mini element to approximate the generalized saddle point problem (\ref{Biot_weakform}). Namely, on discrete level, we have $a_h(\bm{u}_h, \bm{v}_h):=a(\bm{u}_h, \bm{v}_h)$, $b_h(\bm{u}_h, q_h):=b(\bm{u}_h, q_h)$, while the bilinear form for reaction-diffusion operator in pressure space becomes
\begin{equation}\label{d_form_stabilization}
\bar{d}_h(p_h^{n+1},q_h) =d(p_h,q_h) + \epsilon \frac{h^2}{\Delta t} \frac{1}{2\mu+\lambda} \int_{\Omega} \nabla p_h^{n+1} \nabla q_h.
\end{equation}
Here, the second term is the stabilization term. Correspondingly, the right hand side of the second equation is changed to be
$$
(\bar{g}, q_h)=({g}, q_h)+ \epsilon \frac{h^2}{\Delta t} \frac{1}{2\mu+\lambda} \int_{\Omega} \nabla p_h^{n} \nabla q_h.
$$
In the above forms, it is suggested that $\epsilon = \frac{1}{6}$ for inf-sup stable FEs and $\epsilon = \frac{1}{4}$ for equal-order FEs \cite{rodrigo2016stability}. We remark here that the stabilization is quite necessary, because on continuous level if the permeability $\kappa$ is small and due to the Dirichlet boundary condition for the pressure variable, there will be non-physical oscillations in the pressure approximation. Even the conventional inf-sup stable Finite Elements are applied to discretize the Biot model, the resulting linear system does not satisfies the $M$- matrix properties (which means that the discretization scheme is monotone). The stabilization term is applied so that the numerical scheme is monotone \cite{rodrigo2016stability}.

The discrete system obtained from the mixed stabilized finite elements induces the operators $A, B$, and $D$, which are associated with the bilinear forms $a_h(\cdot,\cdot),\ b_h(\cdot,\cdot)$, and $\bar{d}_h(\cdot,\cdot)$ defined in (\ref{a_form1}), (\ref{bilinear_forms_b}), and (\ref{d_form_stabilization}), respectively. The boundedness and the coercivity of $D$ is easy to derive based on the properties of $d(\cdot, \cdot)$ and the stabilization term. We will also denotes $A_{0}$ as the linear operator associated with the vector Laplacian operator in $\bm{V}_h$.

If conforming mixed finite element are applied to discretize the Biot problem, if the the following stability condition \cite{cao2003fast, rodrigo2016stability}
$$
{c}_0{\|q_h\|}^2 \le <( B  A^{-1}  B^t + D)q_h,q_h>  , \quad \forall q_h \in Q_h
$$
holds true, we then have the wellposedness for the discrete problem \cite{cao2003fast, rodrigo2016stability}.
Here, $||\cdot||$ is a properly specified norm, $<\cdot, \cdot>$ is the continuous $L^2$ inner product.

To differentiate the differences of notations, we will also  use ${\bf A}, {\bf B}, {\bf D}$ to denote the resulting matrices. That is,
$$
{\bf A}(i, j) =a_h({\bf b}_i, {\bf b}_j), \quad {\bf B}(i,j)=b_h({\bf b}_i, \phi_j), \quad {\bf D}(i,j)=\bar{d}_h({\phi}_i, \phi_j),  \quad {\bf A}_{0}(i,j)=<\nabla {\bf b}_i, \nabla {\bf b}_j>.
$$
Here, $\{{\bf b}_i\}$ are basis functions of $\bm{V}_h$ and $\{\phi_i\}$ are basis functions of $Q_h$. The properties of the above linear operators and the corresponding matrices can be derived by using the estimates for the bilinear forms. Furthermore, for each vector ${\bf u}$ (or ${\bf p}$), it corresponds to a function $\bm{u}_h$ (or $\bm{p}_h$), and there holds
$$
({\bf A}{\bf u}, {\bf u}) = <A \bm{u}_h, \bm{u}_h>.
$$
Here and thereafter, $(\cdot, \cdot)$ is the discrete $l^2-$ inner product.

\section{Preconditioners and analysis}

In this section, we introduce the preconditioners for the linear system resulted from the stabilized Finite Element discretization of the Biot model. The preconditioners and the spectral properties of the preconditioned system will be presented.

The following block diagonal preconditioner and block triangular preconditioners are frequently used for the saddle point type systems.
\begin{equation}\label{preconditioners}
{\bf P}_1=\left[\begin{array}{cc}
{\bf P}_{\bf A}     & {\bf 0} \\
{\bf 0}   & -{\bf P}_{\bf S}   \\
\end{array}
\right],
\quad \mbox{or} \quad
{\bf P}_2=\left[\begin{array}{cc}
{\bf P}_{\bf A}  & {\bf 0} \\
{\bf B}      & -{\bf P}_{\bf S}   \\
\end{array}
\right],
\quad \mbox{or} \quad
{\bf P}_3=\left[\begin{array}{cc}
{\bf P}_{\bf A}  & {\bf B}^t \\
{\bf 0}          & -{\bf P}_{\bf S}   \\
\end{array}
\right].
\end{equation}
Usually, one can set ${\bf P}_{\bf A}={\bf A}$. The Schur complement for the generalized saddle point system of the form (\ref{stab_sys}) is ${\bf S}={\bf B}{\bf A}^{-1}{\bf B}^t+{\bf D}$.
Based on the Fourier analysis approach \cite{cai2014efficient, cai2015analysis, cai2017fast}, the preconditioner for the Schur complement should be
\begin{equation}\label{Schur_appr}
{\bf P}_{\bf S}= \frac{1}{(2\mu +\lambda)}{\bf M}_p + {\bf D}.
\end{equation}
Here, ${\bf M}_p$, which corresponds to the identity operator in the pressure space, is the pressure mass matrix.
An alternative approximation of the Schur complement is
$$
{\bf P}_{\bf S}= {\bf D},
$$
as ${\bf D}$ is spectral equivalent to ${\bf P}_{\bf S}$ defined in (\ref{Schur_appr}). When combining with Krylov subspace methods, one can use MINRES method for ${\bf M}$ with ${\bf P}_1$, or GMRES method using ${\bf P}_2$ as a left preconditioner or ${\bf P}_3$ as a right preconditioner. We also comment here that if ${\bf P}_{\bf S}={\bf D}$ then the preconditioner is constraint based preconditioner. In the following, we give the estimates based on the exact inverses of ${\bf P}_{\bf A}$ and ${\bf P}_{\bf S}$.  

For the block triangular preconditioner ${\bf P}_2$, it is easy to derive that
\begin{eqnarray}\label{preconditioned_P2M}
{\bf P}_2^{-1} {\bf M} &=&
\left[\begin{array}{cc}
{\bf I}        & {\bf A}^{-1} {\bf B}^t \\
0              & {\bf P}_{\bf S}^{-1} {\bf S}
\end{array}\right].
\end{eqnarray}
Let $\sigma$ and $({\bf u}, {\bf p})^t$ be the eigenvalue and the corresponding eigenvector for (\ref{preconditioned_P2M}), we have
\begin{equation}\label{tri_pre_sys}
\left[\begin{array}{cc}
{\bf A}   & {\bf B}^t \\
{\bf B}   & -{\bf D}
\end{array}\right]
\left[\begin{array}{c}
{\bf u} \\
{\bf p}
\end{array}\right]=
\sigma
\left[\begin{array}{cc}
{\bf A}   & \bm{ 0}  \\
{\bf B}    & -{\bf P}_{\bf S}
\end{array}\right]
\left[\begin{array}{c}
 {\bf u} \\
 {\bf p}
\end{array}\right].
\end{equation}
If $\sigma=1$, we see that ${\bf p}\in \mbox{Ker} ({\bf B}^t)$, and
$$
{\bf B}{\bf u} -{\bf D}{\bf p}={\bf B}{\bf u}-\{ \frac {1}{2\mu+\lambda}~{\bf M}_p +{\bf D} \}{\bf p}.
$$
It follows that ${\bf p}=0$. If $\sigma \neq 1$, there holds ${\bf u}=\frac{1}{\sigma-1}{\bf A}^{-1}{\bf B}^{t} {\bf p}$, plugging into (\ref{tri_pre_sys}), we see that $\sigma$ is the eigenvalue of ${\bf P}_{\bf S}^{-1} {\bf S}$.

As we use conforming finite elements, the matrix properties and the spectral properties of the preconditioned system can be analyzed by estimating the corresponding linear operators. From the estimates of $a(\cdot, \cdot)$, we see that
$$
\alpha_A <A_0 \bm{u}_h, \bm{u}_h> \leq <A \bm{u}_h, \bm{u}_h> \leq (2\mu + \lambda)<A_0 \bm{u}_h, \bm{u}_h>.
$$
Therefore,
\begin{equation}\label{A_quo}
\frac{1}{2\mu+\lambda} ({\bf A}^{-1}_0 {\bf u},{\bf u}) \leq ({\bf A}^{-1} {\bf u}, {\bf u}) \leq \frac{1}{\alpha}_A  ({\bf A}^{-1}_0 {\bf u},{\bf u}), \quad \forall {\bf u}.
\end{equation}
To show that the eigenvalues of ${\bf P}_{\bf S}^{-1} {\bf S}$ have uniform lower and upper bounds independent of mesh refinement and physical parameters, we only need to verify that, $\quad \forall {\bf p}$, the Rayleigh quotient,
\begin{equation}\label{Rayleigh_quo}
\frac{( ({\bf B} {\bf A}^{-1} {\bf B}^{t}+{\bf D}){\bf p},{\bf p})}{\left( (\frac{1}{{2\mu+\lambda}}{\bf M}_{p} +{\bf D}) {\bf p}, {\bf p} \right)}
\end{equation}
has uniform lower and upper bounds. From (\ref{A_quo}), if ${\bf u}={\bf B}^t {\bf p}$, we have
$$
\frac{1}{2\mu+\lambda} ({\bf B} {\bf A}_0^{-1}{\bf B}^t {\bf p}, {\bf p}) \le ({\bf B}{\bf A}^{-1} {\bf B}^t {\bf p}, {\bf p}) \le \frac{1}{\alpha_A} ({\bf B} {\bf A}_0^{-1} {\bf B}^t {\bf p}, {\bf p}).
$$
Plugged into (\ref{Rayleigh_quo}), we have
\begin{equation}\label{est_Schur}
\frac{\frac{1}{2\mu+\lambda} ({\bf B} {\bf A}_0^{-1}{\bf B}^t {\bf p}, {\bf p})+({\bf D}{\bf p}, {\bf p})}{\left( (\frac{1}{{2\mu+\lambda}}{\bf M}_{p} +{\bf D}) {\bf p}, {\bf p} \right)} \le \frac{( ({\bf B} {\bf A}^{-1} {\bf B}^{t}+{\bf D}){\bf p},{\bf p})}{\left( (\frac{1}{{2\mu+\lambda}}{\bf M}_{p} +{\bf D}) {\bf p}, {\bf p} \right)} \le \frac{\frac{1}{\alpha_A} ({\bf B} {\bf A}_0^{-1} {\bf B}^t {\bf p}, {\bf p})+ ({\bf D}{\bf p}, {\bf p})}{\left( (\frac{1}{{2\mu+\lambda}}{\bf M}_{p} +{\bf D}) {\bf p}, {\bf p} \right)}.
\end{equation}
On the other hand, as we use inf-sup stable conforming Finite Elements, there exists a $\beta >0$, independent of mesh refinement, such that
$$
\displaystyle
\beta \le  \inf_{p_h \in Q_h} \sup_{\bm{v}_h \in \bm{V}_h} \frac{b(\bm{v}_h, p_h)}{||\bm{v}_h||_{H^1}||p_h||_{L^2}}.
$$
Combined with the facts that the operators $A_0$ and $B$ are bounded, we have the following spectral properties of the matrices \cite{elman2005finite, cai2015analysis}:
\begin{equation}\label{schur_est_Lap}
\beta^2 \le \frac{({\bf B} {\bf A}_0^{-1} {\bf B}^t {\bf p}, {\bf p})}{({\bf M}_p {\bf p}, {\bf p})} \le 1.
\end{equation}
Plugging the inequality (\ref{schur_est_Lap}) into (\ref{est_Schur}), we see clear that the eigenvalues of ${\bf P}_{\bf S}^{-1} {\bf S}$ have uniform lower and upper bounds independent of mesh refinement and physical parameters. In summary, we have the following proposition.

\begin{theorem}
For the Biot problem discretized by the stabilized Mini elements, if exact elliptic solvers are applied in ${\bf P}_{2}^{-1}$ (or ${\bf P}_3^{-1}$), the eigenvalues of the preconditioned system are either $1$ (with the multiplicity equals to the number of unknowns for the displacement variables) or have uniform lower and upper bounds independent of mesh refinement and physical parameters.
\end{theorem}

{\bf Remark}. We remark here that no matter there is a stabilization term or not, ${\bf P}_{2}^{-1}{\bf M}$ has eigenvalues which have uniform lower and upper bounds independent of mesh refinement and physical parameters. However, it is crucial to add the stabilization term to ensure the monotonicity of the discretization especially when the permeability is small.

%

\section{Iterative algorithms}

\subsection{A preconditioned GMRES method}
For saddle point problems and generalized saddle point problems, a commonly used solution strategy is applying the GMRES algorithm as the outer iteration, a preconditioner, implemented using Multigrid solvers \cite{xu1992iterative} or domain decomposition solvers \cite{toselli2005domain, dohrmann2009overlapping, dohrmann2010hybrid} or spare solvers, is employed as the inner iteration. The detailed algorithm of GMRES method can be found in \cite{saad1986gmres, saad1993flexible}. In this paper, we also use GMRES method as the outer iteration method. The preconditioner takes the form of (\ref{preconditioners}) with ${\bf P}_{\bf A}^{-1}$ and ${\bf P}_{\bf S}^{-1}$ being implemented as incomplete Cholesky factorizations. More efficient and advanced implementation can be based on overlapping Schwarz domain decomposition methods \cite{cai2015overlapping, cai2016hybrid} or Multigrid method \cite{gaspar2007efficient, cai2014efficient, cai2015analysis} for the original or the mixed reformulation of the linear elasticity operator. In our implementation, both (almost) exact and inexact solves of ${\bf P}_{\bf A}^{-1}$ and ${\bf P}_{\bf S}^{-1}$ are applied to check the effects of inexact solve, see the numerical experiments in Section 5 for more details.

\subsection{Variable-relaxation parameter Uzawa algorithm.}

It is well-known that the classical Uzawa algorithm \cite{bacuta2006unified, benzi2005numerical, bramble1997analysis, uzawa1958iterative} converges slowly, and is not of practical use in many applications. A lot of works have been done to accelerate the efficiency, see \cite{elman1994inexact, hu2001iterative, hu2006nonlinear, hu2002two, ho2015accelerating} and the references cited therein. For example, the authors of \cite{hu2001iterative, hu2006nonlinear, hu2002two} consider the classical saddle point problem (with the $(2, 2)$ block being zero), they introduced some variable relaxation parameters in each step of Uzawa algorithm so that the errors are minimized under proper energy norms in each step. It is natural to apply the variable-relaxation parameter Uzawa algorithm to the generalize saddle point system studied in this paper. The variable-relaxation parameter Uzawa algorithm for the generalized saddle point system reads as:
\begin{equation}\label{Uzawa_Alg1}
 \left \{
 \begin{array} {ll}
	{\bf u}_{k+1}={\bf u}_{k}+ \omega_k {\bf P}_{\bf A}^{-1} [{\bf f}- ({\bf A} {\bf u}_k+{\bf B}^{t} {\bf p}_k)], \\
    {\bf p}_{k+1}={\bf p}_{k}- \theta_k \tau_k   {\bf P}_{\bf S}^{-1} [{\bf g}-({\bf B}{\bf u}_{k+1}-{\bf D}{\bf p}_k)].
	\end{array}
    \right.
\end{equation}
The parameter $\omega_k$ is chosen so that the error ${\bf u}-{\bf u}_k$ is minimized under the ${\bf A}$- norm. Denoting
${\bf f}_k = {\bf f} -({\bf A} {\bf u}_k+ {\bf B}^t {\bf p}_k), {\bf c}_k = {\bf P}_{\bf A}^{-1}{\bf f}_k$, then a prototype
choice of $\omega_k$ is
$$
\omega_k= \frac{({\bf f}_k,{\bf c}_k)}{({\bf A} {\bf c}_k,{\bf c}_k)}.
$$
The parameter $\tau_k$ is chosen so that the error ${\bf p}-{\bf p}_k$ is minimized under the ${\bf P}_{\bf S}$- norm or its equivalent norm. Denoting ${\bf g}_k = {\bf g}-({\bf B} {\bf u}_{k+1}-{\bf D}{\bf p}_k) , {\bf d}_k= -{\bf P}_{\bf S}^{-1} {\bf g}_k$, then a prototype choice of $\tau_k$ is
$$
\tau_k = \frac {({\bf g}_k,{\bf d}_k)}{({\bf P}_{\bf S} {\bf d}_k, {\bf d}_k)} \quad \mbox{or} \quad \tau_k= \frac {({\bf g}_k,{\bf d}_k)}{({\bf S} {\bf d}_k, {\bf d}_k)}.
$$
In the denominator part, ${\bf P}_{\bf S}$ can be replaced by any spectral equivalent matrix. For the parameter $\theta_k$, it is a damping parameter so that the convergence of the algorithm can be guaranteed \cite{hu2001iterative}. If there is no $(2,2)$ block, it is suggested in \cite{hu2001iterative} that
$$
{\theta}_k = \frac {1- \sqrt{1- {\omega}_k}}{2}.
$$
However, from our numerical experience, $\theta_k=1.0$ usually leads to the best performance if almost exact Poisson solvers are used for ${\bf P}_{\bf A}$ and ${\bf P}_{\bf S}$, no matter whether there is a $(2,2)$ block. In summary, the algorithm for variable relaxation parameter Uzawa algorithm is as listed in {\bf Algorithm 1}.

\begin{algorithm}\label{Alg_UA}
\caption{Uzawa algorithm with variable relaxation parameters.}\label{multiplicative_alg}
\begin{algorithmic}[1]
\State Given the initial guesses ${\bf u}_0 \in {R}^n$ and ${\bf p}_0 \in
 {R}^m$ compute the sequences {${\bf u}_k$,${\bf p}_k$} for i = 1,2,... as follows.
\State Step 1. Compute ${\bf f}_k = {\bf f} -({\bf A} {\bf u}_k+ {\bf B}^t {\bf p}_k), {\bf c}_k = {\bf P}_{\bf A}^{-1}{\bf f}_k$, and
$$
{\omega}_k=\left\{
\begin{array}{lll}  
\frac {({\bf f}_k,{\bf c}_k)}{({\bf A} {\bf c}_k,{\bf c}_k)}, \quad {\bf f}_k \neq 0, \\
1,~~~~~~~~~ \quad {\bf f}_k = 0.\\
\end{array}
\right.
$$
\State Set ${\bf u}_{k+1} = {\bf u}_k + {\omega}_k {\bf c}_k$

\State Step 2. Compute ${\bf g}_k = {\bf g}-({\bf B} {\bf u}_{k+1}-{\bf D}{\bf p}_k) , {\bf d}_k= -{\bf P}_{\bf S}^{-1} {\bf g}_k$ and,
$$
{\tau}_k=\left\{
\begin{array}{lll}
\frac {({\bf g}_k,{\bf d}_k)}{({\bf P}_{\bf S} {\bf d}_k, {\bf d}_k)}, \quad {\bf g}_k \neq 0, \\
1,       ~~~~~~~~~~~ \quad {\bf g}_k = 0.\\
\end{array}
\right.
$$
\State Set
${\bf p}_{k+1} = {\bf p}_k + {\tau}_k{\bf d}_k$.

\end{algorithmic}
\end{algorithm}

We further comment here that if $\omega_k=1.0$ and $\theta_k \tau_k=1.0$ in (\ref{Uzawa_Alg1}), the above algorithm is the classical preconditioned Uzawa algorithm. Similar to the classical Uzawa algorithm, the performance of these algorithms depends on the estimates of the maximal and minimal eigenvalues of the preconditioned Schur complement \cite{hu2001iterative, hu2006nonlinear, hu2002two}.

\subsection{Anderson accelerated Uzawa algorithm.}

There have been some successful applications of Anderson acceleration algorithm for solving nonlinear problems \cite{walker2011anderson}. Some recent progress on Anderson acceleration algorithm can be found in \cite{walker2011anderson, toth2015convergence, ho2015accelerating}. In a theoretical paper \cite{toth2015convergence}, it is shown that Anderson acceleration algorithm is equivalent to the GMRES method in a certain sense when they are applied to a single linear system. In this part of research, we intend to generalize and apply the method to Biot's model. However, there is not enough careful comparison to clarify whether the Anderson acceleration algorithm is superior to the GMRES method or other iterative methods. In \cite{ho2015accelerating}, the authors show that the preconditioned accelerated Uzawa algorithm is comparable to other algorithms for Oseen equations if exact Poisson solvers are applied in preconditioning steps. However, the comparisons of these algorithms based on inexact Poisson solvers are not provided, even for the linear Stokes problem.

For a general fixed-point iteration,
$$
x_{k+1}=g(x_k),
$$
By employing the results of previous steps, the Anderson acceleration algorithm can provide a better approximation of the true solution \cite{walker2011anderson, toth2015convergence}. The Anderson acceleration is described in {\bf Algorithm 2}. There have been some successful applications of Anderson acceleration in many nonlinear problems, see \cite{walker2011anderson} and the references cited therein. We comment here that Anderson acceleration can be applied to any fixed point iteration, no matter the iteration is linea or nonlinear. In this part of the project, we intend to apply Uzawa algorithm and treat it as a fixed point iteration for \eqref{Uzawa_fixed}. Anderson acceleration is employed to speed up the convergence rate of the fixed point iteration. Preconditioners will be used in each step of the Uzawa algorithm to improve the efficiency.

\begin{algorithm}\label{Alg_AA}
\caption{Anderson acceleration (AA) of fixed point iteration.}\label{multiplicative_alg}
\begin{algorithmic}[1]
\State Given ${x}_0$ and $m \ge 1$. Set $x_1 = g(x_0)$.
\State For $k = 1, 2, . . .$ \par
\State ~~~~~Set $m_k = \min \{m, k\}$.
\State ~~~~~Set $F_k = (f_{k-m_k}, . . . , f_k)$, where $f_i = g(x_i)-x_i$.
\State ~~~~~Determine $\alpha^{(k)} = (\alpha^{(k)}_0, . . ., \alpha^{(k)}_{m_k} )^t$ that solves
$$
\min_{\alpha =(\alpha_0, . . ., \alpha_{m_k} )^t} ||F_k \alpha ||_2 ~~ \mbox{s.t.}~~ \sum_{i=0}^{m_k} \alpha_i=1.
$$
\State ~~~~~Set $x_{k+1} = \sum_{i=0}^{m_k} \alpha^{(k)}_i g(x_{k-m_k+i})$.
\end{algorithmic}
\end{algorithm}

In a recent work \cite{ho2015accelerating}, the authors propose to apply Anderson acceleration \cite{walker2011anderson, toth2015convergence} to improve the performance of Uzawa algorithm for saddle point problems. To apply the Anderson acceleration to the generalized saddle point system studied in this work, we can firstly rewrite the Uzawa algorithm as a fixed point iteration as follows.
\begin{equation}
\left\{
\begin{array}{lll}  
{\bf u}_{k+1} &=&{\bf u}_{k}  + {\bf P}_{\bf A}^{-1} ({\bf f} - {\bf A} {\bf u}_{k} -{\bf B}^t {\bf p}_{k}), \\
{\bf p}_{k+1}   & =  &{\bf p}_{k} +{\bf P}_{\bf S}^{-1} ({\bf g}-{\bf B}{\bf u}_{k+1} +{\bf D} {\bf p}_{k}). \\
\end{array}\right.
\label{Uzawa_fixed}
\end{equation}
Here, ${\bf P}_{\bf A}^{-1}$ and ${\bf P}_{\bf S}^{-1}$ may involve relaxation parameters or scaling factors. Rewriting (\ref{Uzawa_fixed}) as a fixed point iteration, we obtain
$$
\left[\begin{array}{cc}
{\bf P}_{\bf A}    & {\bf 0} \\
{\bf B}            & {\bf P}_{\bf S}   \\
\end{array}
\right]
\left[\begin{array}{c}
{\bf u}_{k+1}  \\
{\bf p}_{k+1}       \\
\end{array}
\right] =
\left[\begin{array}{cc}
{\bf P}_{\bf A}{\bf -A}   & -{\bf B}^t \\
{\bf 0}         & {\bf P}_{\bf S}-{\bf D}   \\
\end{array}
\right]
\left[\begin{array}{c}
{\bf u}_{k}     \\
{\bf p}_{k}      \\
\end{array}
\right]+
\left[\begin{array}{c}
{\bf f}     \\
{\bf g}      \\
\end{array}
\right].
$$
By employing the results of previous steps, we expect that the Anderson acceleration algorithm
can provide a better approximation of the true solution than the typical Uzawa algorithm itself. We will also check whether the combination of Anderson acceleration with Uzawa leads to good performance for the generalized saddle point problem studied in this paper.

\section{Numerical experiments}

We compare all the algorithms in this section. The computational domain is $[0, 1]^2$, we vary the meshsize and the physical parameters, so that we can test the robustness of our preconditioners and the effects of mesh refinement, physical parameters. In our tests, we set the Poisson ratio $\nu=0.3$ for testing the compressible case and set $\nu=0.49$ for testing the almost incompressible case. The Young's module $E$ is fixed to be $1000$. The permeability is set to be $1$ for testing the large permeability case and is set to be $0.0001$ for testing the small permeability case. Numerical experiments are summarized in Table \ref{Compare_Algs1} to Table \ref{Compare_Algs3}. In our tests, the stopping tolerances are either the relative errors or the relative residuals in $l^2$ norm need to be reduced to $1.0\times 10^{-6}$. We record the number of iterations used for different algorithms. "DOFs" means the total number of degrees of freedom, "Nx" means the number of elements along each direction, "No Pre" means GMRES algorithm with no preconditioner being used, "PGMRES" means the preconditioned GMRES method, "U" means the preconditioned Uzawa with the exact Poisson solvers, "AAU" means the Anderson accelerated Uzawa algorithm, "VU" means the variable relaxation parameter Uzawa algorithm. Correspondingly, "IPGMRES" means the preconditioned GMRES method with inexact Poisson solvers, "IU" means the preconditioned Uzawa algorithm with inexact Poisson solvers, "IAAU" and "IVU" means the Anderson accelerated Uzawa algorithm and the variable relaxation parameter Uzawa algorithm with inexact Poisson solvers respectively. For Uzawa algorithm and the inexact Uzawa algorithm, the relaxation parameter $\omega=2.5$ is an empirical choice (our extensive experiments show that $\omega=2.5$ is almost the optimal) .

First of all, from Table \ref{Compare_Algs1} to Table \ref{Compare_Algs3}, we see clearly that if the Poisson solvers are implemented exactly, it only needs several iterations (GMRES only needs $4$ to $5$ iterations, all the other algorithms also need less than $10$ iterations). These results clearly show that our preconditioners are very effective for the Biot operator. Moreover, from Table \ref{Compare_Algs1} to Table \ref{Compare_Algs3}, these preconditioners are very robust with respect to the mesh refinement and the physical parameters (not only the Poisson ratio, but also the permeability). Moreover, from the results, Anderson accelerated Uzawa (AAU) gives better performance than the Uzawa (U) algorithm and the variable relaxation Uzawa (VU) algorithm. By comparing with the results obtained by using the inexact Poisson solvers, it is obvious that the inexact solvers will make the total number of iterations much larger. We will explain more details on the effects of inexact Poisson solvers in the following.

In our implementation, although the matrix ${\bf D}$ is ill conditioned, the incomplete Cholesky factorization is not very difficult to apply because its condition number mainly depends on mesh refinement and will not be affected too much by the physical parameter. In our implementation,
$$
{\bf D} = {\bf L}_{\bf D} {\bf L}^t_{\bf D} + tol
$$
with the tolerance $tol=0.001$. However, for the matrix ${\bf A}$, its condition number not only depends on mesh refinement but also depend on the physical parameter, in particular, $\lambda$. When the Poisson ratio $\nu$ approaches to $0.5$, $\lambda$ becomes very large, and the condition number of ${\bf A}$ is huge. For such kind of matrix, incomplete Cholesky factorization does not work very well (although the matrix itself is symmetric positive definite). To make the inexact solvers works, we apply the modified incomplete Cholesky factorization: applying the incomplete factorization to
$$
{\bf  A} + \alpha * \mbox{diag}(\mbox{diag}({\bf A})) = {\bf L}_{\bf A} {\bf L}^t_{\bf A} + tol.
$$
In our implementation, we set $\alpha=10$ and the tolerance is also set to be $0.001$ (cf. Matlab function "ichol.m"). As an example to show that inexact Poisson solvers with Anderson accelerated Uzawa (AA Uzawa) makes different, we check the preliminary results in Table \ref{Compare_Algs1} to Table \ref{Compare_Algs3}. If inexact Poisson solvers (${\bf A}^{-1}$ and ${\bf P}_{\bf S}^{-1}$ are approximated by using their incomplete Cholesky factorizations or modified incomplete Cholesky factorizations) are applied, the advantages of Anderson accelerated Uzawa algorithm with inexact Poisson solvers (IAAU) is more obvious when compared with the Uzawa algorithm and the varaible-relaxation parameter Uzawa algorithm with the inexact Poisson solvers. For both the Uzawa algorithm and the variable-relaxation parameter Uzawa algorithm, when inexact Poisson solvers are applied, their performance are not good. Although the variable-relaxation parameter Uzawa algorithm gives slightly better performance than that of Uzawa algorithm, its performance is much worse than that of Anderson accelerated algorithm or GMRES algorithm.

\begin{table}[h]
\begin{center}
  \centering
  \begin{tabular}{r|r||ccccc||cccc}
  \hline
{DOFs} &  Nx & No Pre & PGMRES &U  &AAU &VU  &IPGMRES & IU    &IAAU    &IVU \\
 \hline
  $1891$  &16  &475        &4    &8   &5     &8   &63   &$>3000$  &1629   &762  \\
 $7363$   &32  &1462       &4    &8   &6     &8   &116  &$>3000$ &$>3000$ &2880  \\
 $29059$  &64  &$>3000$    &4    &8   &6     &8   &222  &$>3000$ &$>3000$ &$>3000$ \\
 $115459$ &128 &$>3000$    &4    &8   &6     &9   &450  &$>3000$ &$>3000$ &$>3000$ \\
 \hline
  \end{tabular}
\caption{Numbers of iterations using different algorithms with exact (U, AAU, and VU) Poisson solvers and inexact Poisson solvers (IU, IAAU, and IVU). $E=1000$, $\nu=0.3$, and all the other parameters are equal to $1$.
}\label{Compare_Algs1}
\end{center}
\end{table}

\begin{table}[h]
\begin{center}
  \centering
  \begin{tabular}{r|r||ccccc||cccc}
  \hline
{DOFs} &  Nx & No Pre & PGMRES &U  &AAU &VU  &IPGMRES & IU    &IAAU    &IVU \\
 \hline
  $1891$  &16  &533     &4    &8   &5     &8   &111   &$>3000$   &1922    &2334  \\
  $7363$  &32  &$>3000$ &4    &8   &5     &8   &197   &$>3000$   &$>3000$    &$>3000$ \\
 $29059$  &64  &$>3000$ &4    &8   &5     &8   &321   &$>3000$   &$>3000$    &$>3000$ \\
$115459$ &128  &$>3000$ &4    &9   &5     &8   &515    &$>3000$   &$>3000$    &$>3000$\\
 \hline
  \end{tabular}
\caption{Numbers of iterations using different algorithms with exact (U, AAU, and VU) Poisson solvers and inexact Poisson solvers (IU, IAAU, and IVU). $E=1000$, $\nu=0.49$ and all the other parameters are equal to $1$.
}\label{Compare_Algs2}
\end{center}
\end{table}

\begin{table}[h]
\begin{center}
  \centering
  \begin{tabular}{r|r||ccccc||cccc}
  \hline
{DOFs} &  Nx & No Pre & PGMRES &U  &AAU &VU  &IPGMRES & IU    &IAAU    &IVU \\
 \hline
  $1891$  &16  &269      &4    &8    &5    &8    &53   &$>3000$  &935    &859  \\
 $7363$   &32  &543      &4    &8    &6    &8    &78   &$>3000$  &726   &$>3000$\\
 $29059$  &64  &$>3000$  &4    &8    &6    &8    &139  &$>3000$  &1332   &$>3000$ \\
 $115459$ &128 &$>3000$  &4    &9    &6    &9    &221  &$>3000$  &2336   &$>3000$ \\
 \hline
  \end{tabular}
\caption{Numbers of iterations using different algorithms with exact (U, AAU, and VU) Poisson solvers and inexact Poisson solvers (IU, IAAU, and IVU). $E=1000$, $\nu=0.3$, $\kappa=0.0001$ and all the other parameters are equal to $1$.
}\label{Compare_Algs3}
\end{center}
\end{table}

By comparing the results obtained by all the algorithms, we also observe that the performance of Anderson acceleration Uzawa algorithm is not that good when inexact Poisson solvers are used. If exact Poisson solvers are employed in preconditioners, the performance of Anderson acceleration Uzawa algorithm is comparable with that of GMRES method. However, when the inexact Poisson solvers are employed, the performance of Anderson acceleration Uzawa algorithm is much worse than that of GMRES method. Another observation is that the performance of variable-relaxation parameter Uzawa algorithm is not good either. Actually, this algorithm is somehow a kind of Krylov subspace method (as the errors are minimized in each sub-step). However, the variable-relaxation parameter Uzawa algorithm needs many iterations to converge especially when the Poisson ratio is close to $0.5$ (cf. Table \ref{Compare_Algs2}).

When inexact Poisson solvers are employed, the physical parameters do have critical effects on the performance of the iterative algorithms. By comparing the results from Table \ref{Compare_Algs1} and Table \ref{Compare_Algs2}, we see that when the elastic material is almost incompressible, more iterations are needed for all iterative algorithms compared with those for compressible elastic material. By comparing the results from Table \ref{Compare_Algs1} and Table \ref{Compare_Algs3}, we see that small permeability actually does not make too much trouble to the performance of the all algorithms. In fact, when the permeability becomes smaller, the numbers of iterations needed are less.

\section{Conclusions}

In this paper, we intend to develop efficient iterative methods and robust preconditioners for solving Biot equations. Our investigation is very comprehensive and in details. Both exact Poisson solvers and inexact Poisson solvers are employed in in preconditioning steps. We compare the GMRES method, Uzawa method, Anderson accelerated Uzawa algorithm and also the variable-relaxation parameter Uzawa algorithm \cite{hu2001iterative, hu2006nonlinear, hu2002two}. We conduct numerical analysis and experiments to highlight the advantages and the disadvantages of each algorithm.

From the numerical experiments, it is observed that the GMRES method combined with block triangular preconditioner still gives the best convergence rate. The possible reason is that GMRES method has the Galerkin property and minimize the {\bf global} residual in each step of the iteration. No matter exact or inexact Poisson solvers are used, it needs the least number of iterations among all the algorithms. From our investigation, the advantage of combining Anderson acceleration with Uzawa algorithms is not obvious. We predict that even for Stokes problem or Oseen problem, there is no obvious benefit of combining Anderson acceleration with Uzawa algorithm over the GMRES method using block triangular preconditioners especially when inexact Poisson solvers are employed. For variable relaxation Uzawa algorithm, it is actually also a Krylov subspace method. Although the convergence properties seem to be not as good as the GMRES method for saddle point problem, the method is a memory saving method. In contract, both GMRES method and Anderson acceleration algorithm are memory cost approaches. For Anderson acceleration, although we do not have very positive conclusion for the Biot model studied in this paper, it is still very promising to apply the algorithm to nonlinear problems, in particular, those nonlinear fixed-point problems with discontinuous Jacobians.

\renewcommand{\refname}{References}

\end{document}